\documentclass[12pt]{amsart}
\usepackage{graphicx}
\usepackage[frenchb]{babel}
\usepackage{enumerate}
\input xy
\xyoption{all}
\usepackage[latin 1]{inputenc}
\newtheorem{thm}{Th\'{e}or\`{e}me}[section]
\newtheorem{cor}[thm]{Corollaire}

\newtheorem{prob}[thm]{Probl\`{e}me}
\newtheorem{conj}[thm]{Conjecture}
\newtheorem{lem}[thm]{Lemme}
\newtheorem{prop}[thm]{Proposition}
\newtheorem{defn}[thm]{D\'{e}finition}
\theoremstyle{remark}

\newtheorem{rem}[thm]{Remarque}
\numberwithin{equation}{section}
\begin{document}

\title{Sur la cat\'egorie des bimodules de Soergel}%
\author{Nicolas Libedinsky}%
\dedicatory{Para Javi, mi hoja ligera}

\address{UFR de Math\'ematiques et Institut de Math\'ematiques de
 Jussieu, Universit\'e Paris 7, 2 place Jussieu, 75251 Paris Cedex 05,
 France}%
\email{libedinsky@math.jussieu.fr}%

\begin{abstract}
 La cat\'egorie $\mathbf{B}$ de Soergel  d'un syst\`eme de Coxeter
 $(W,\mathcal{S})$ est une cat\'egorie de bimodules sur une alg\`ebre de polynômes sur
 laquelle $W$ agit. C'est une cat\'egorification   de l'alg\`ebre de Hecke de
 $(W,\mathcal{S})$. Dans cet article nous donnons une description
 combinatoire des espaces de morphismes dans $\mathbf{B}$. En corollaire, on obtient une description analogue des morphismes dans $\mathcal{O}_0$-proj, où $\mathcal{O}_0$ est le bloc principal de la catégorie $\mathcal{O}$ de \textbf{BGG}.
\end{abstract}

\maketitle

\section{Introduction}
 En 1980, Kazhdan et Lusztig ont pos\'e leur c\'el\`ebre conjecture de
 positivit\'e \cite{KL}. Si $(W,\mathcal{S})$ est un syst\`eme de
 Coxeter et $\mathcal{H}$ son alg\`ebre de Hecke, cette conjecture dit que
 certains polyn\^omes de changement de base dans $\mathcal{H}$ (les
 polyn\^omes de Kazhdan-Lusztig) ont des coefficients positifs. Ces polyn\^omes
 ont donn\'e naissance \`a ce qu'on appelle la th\'eorie de
 Kazhdan-Lusztig, qui s'est av\'er\'ee avoir des  liens profonds avec la
 g\'eom\'etrie, les repr\'esentations et la combinatoire.

 En 1992 \cite{S1}
 Soergel a cat\'egorifi\'e $\mathcal{H}$, c'est \`a dire qu'il a d\'efini
 une cat\'egorie tensorielle $\mathbf{B}$ et un isomorphisme d'anneaux $\varepsilon$
 de $\mathcal{H}$ vers le groupe de Grothendieck scind\'e de
  $\mathbf{B}$. Il a alors pos\'e une conjecture
 (\ref{cs} ci-dessous) qui relie, via $\varepsilon$, les \'el\'ements de la
 base de Kazhdan-Lusztig avec les \'el\'ements ind\'ecomposables de
 $\mathbf{B}$. Cette conjecture implique la conjecture de positivit\'e de
 Kazhdan-Lusztig.

 Historiquement, Soergel a introduit cette cat\'egorie dans le cas des
 groupes de Weyl pour relier la cat\'egorie $\mathcal{O}$ d'une alg\`ebre de
 Lie semisimple complexe et les faisceaux pervers sur les vari\'et\'es de
 drapeaux. Le lien avec la g\'eom\'etrie lui a permis de d\'emontrer sa conjecture dans le cas des groupes de Weyl, et il en a d\'eduit une nouvelle preuve de la conjecture de
 Kazhdan-Lusztig, d\'emontr\'ee auparavant par Beilinson et Bernstein
 \cite{BB} et par Brylinski et Kashiwara \cite{BK}.

 Par ailleurs, dans un article r\'ecent \cite{K}, Khovanov montre qu'\`a partir de certains complexes dont les termes appartiennent \`a $\mathbf{B}$, on retrouve  l'homologie r\'eduite d\'efinie dans \cite{KR}, qui est un invariant de noeuds et dont la caract\'eristique d'Euler  est le polyn\^ome HOMFLYPT.

Dans la section 2 de cet article, nous donnons l'\'enonc\'e pr\'ecis du
 th\'eor\`eme de Soergel, de la conjecture de positivit\'e de
 Kazhdan-Lusztig ainsi que de la conjecture de Soergel, en expliquant \`a chaque
 fois quelles parties de ces conjectures ont \'et\'e prouv\'ees. Dans
  tous les cas (sauf pour les groupes de Coxeter universels), les parties de ces
 conjectures qui sont r\'esolues, sont les cas o\`u on peut associer une
 g\'eom\'etrie. C'est pour cela en partie que cette approche alg\'ebrique
 de Soergel est particuli\`erement int\'eressante.

Dans la section 3, nous exposons les r\'esultats d\'ej\`a connus qui
 nous int\'eressent sur les espaces d'homomorphismes de la cat\'egorie
 $\mathbf{B}$. Tous ces r\'esultats sont explicites ou implicites dans
 l'article de Soergel \cite{S3}. Le th\'eor\`eme fondamental de cette
 section est le th\'eor\`eme \ref{inversa}, qui donne les dimensions gradu\'ees de
 ces espaces d'homomorphismes (corollaire \ref{grados}).

 Dans la section 4, nous d\'efinissons la
 \og base des feuilles l\'eg\`eres \fg\ (BFL), qui est un sous-ensemble de ces
 espaces d'homomorphismes.   On peut regarder la construction de la section 4 comme une cat\'egorification de la formule donn\'ee dans
 le corollaire \ref{grados}. Cette formule dit que les dimensions gradu\'ees des espaces d'homomorphismes sont donn\'ees par le coefficient en $1$ d'un produit du type $(1+T_{s_1})\cdots (1+T_{s_n})$. 
 
 Dans les sous-sections 4.1-4.3 nous d\'efinissons un morphisme entre certains \'el\'ements de
 $\textbf{B}$ associ\'e \`a une relation de tresses. Dans les sous-sections 4.4 et 4.5 nous  construisons par r\'ecurrence la BFL, en imitant au niveau des morphismes de la cat\'egorie $\textbf{B}$ la r\'ecurrence qui appara\^{\i}t pour le calcul du produit $(1+T_{s_1})\cdots (1+T_{s_n})$ \`a partir du produit $(1+T_{s_1})\cdots (1+T_{s_{n-1}})$.

 Dans la section 5, nous prouvons le th\'eor\`eme fondamental de cet article : la BFL est en fait, comme le laisse pr\'evoir son nom, une base de l'espace
 d'homomorphismes, comme module \`a droite sur un certain anneau.  Dans la section 6 nous en d\'eduisons des bases pour les morphismes dans $\textbf{B}$. 

Soit $\mathcal{O}_0$ le bloc principal de la catégorie $\mathcal{O}$ de \textbf{BGG}. Comme corollaire du résultat de la section 6,  nous trouvons explicitement,  dans la section 7, les morphismes  dans la sous-catégorie pleine de $\mathcal{O}_0$ d'objets projectifs ($\mathcal{O}_0-$proj).
 
La catégorie $\mathbf{B}$ est construite à partir d'une représentation réflexion fidèle de $W$. Dans l'article en préparation \cite{lib} nous utilisons le théorème \ref{bacan} pour prouver des équivalences entre les conjectures de Soergel  associées à  différentes représentations de $W$. En particulier nous montrons que sur $\mathbb{R}$ il suffit de considérer la représentation géométrique. 

\hspace{20mm}

Je tiens \`a remercier Raphaël Rouquier pour son aide dans la r\'edaction de ce travail.

\section{Conjecture de Soergel}

\subsection{Th\'{e}or\`{e}me de Soergel}
Le th\'{e}or\`{e}me de Soergel affirme que la cat\'egorie de Soergel
 est une cat\'egorification de l'alg\`ebre de Hecke d'un groupe de
 Coxeter. Donnons quelques d\'{e}finitions avant de donner l'\'enonc\'e
 pr\'ecis de ce th\'{e}or\`{e}me.

\begin{defn}
Un syst\`{e}me de Coxeter est un couple $(W,\mathcal{S})$ o\`u $W$
est un groupe et $\mathcal{S}\subseteq W$ une partie
g\'{e}n\'{e}ratrice, tels que $W$ admet une pr\'{e}sentation de
g\'{e}n\'{e}rateurs $s\in \mathcal{S}$ et relations
$(sr)^{m(s,r)}=1$ pour $s,r\in \mathcal{S}$, avec $m(s,s)=1$ ,
$m(s,r)\geq 2$ et \'{e}ventuellement $m(r,s)=\infty$ si $s\neq r.$
\end{defn}

\begin{defn}
Soit  $(W, \mathcal{S})$ un syst\`{e}me de Coxeter. On d\'{e}finit
l'alg\`{e}bre de Hecke $\mathcal{H} = \mathcal{H}(W, \mathcal{S})$
comme la $\mathbb{Z} [v,v^{-1}]$-alg\`ebre de g\'{e}n\'{e}rateurs $\{
 T_{s}
\}_{s\in \mathcal{S}}$, ceux-ci satisfaisant les relations $$T^{2}_{s}=
v^{-2}+(v^{-2}-1)T_{s}$$ pour tout $s\in \mathcal{S}$ et
$$\underbrace{T_{s}T_{r}T_{s}...}_{m(s,r)\, \mathrm{termes}
}=\underbrace{T_{r}T_{s}T_{r}...}_{m(s,r)\, \mathrm{termes} }$$ si
$s,r \in \mathcal{S}$ et $sr$ est d'ordre $m(s,r)$.

Si
$x=s_{1}s_{2}\cdots s_{n}$ est une expression r\'{e}duite de $x$, on
d\'{e}finit $T_x=T_{s_1}T_{s_2}\cdots T_{s_n}$ ($T_x$ ne d\'epend pas du choix de la d\'ecomposition r\'eduite). On pose $q=v^{-2}.$

\end{defn}

Soit $\mathcal{T} \subseteq W$ le sous-ensemble des \og r\'eflexions \fg,
c'est \`{a} dire, tous les \'{e}l\'{e}ments qui sont conjugu\'{e}s
aux \'{e}l\'{e}ments de $\mathcal{S}$.

\begin{defn}
Une repr\'esentation de dimension finie de $W$ (sur un corps $k$ de caract\'eristique $\mathrm{car}$(k) diff\'erente de $2$) est
appel\'{e}e r\'eflexion fid\`{e}le (RF) si elle est fid\`{e}le et
 si les  \'{e}l\'{e}ments de $W$ qui ont un espace de points fixes de
codimension un forment exactement l'ensemble de r\'eflexions.
\end{defn}

 \begin{rem}Dans \cite{S3}  Soergel montre que si $k=\mathbb{R}$, pour  tout $W$ il existe au moins une RF.
 \end{rem}
  Soit $V$ une RF. Soit
$R=S(V^{*})= R(V)$ l'alg\`{e}bre sym\'{e}trique de $V^*$, c'est \`a
 dire
l'alg\`{e}bre des fonctions r\'{e}gulieres sur $V$, sur laquelle $W$
agit par fonctorialit\'{e}.

\begin{defn}
Pour toute petite cat\'{e}gorie additive $\mathcal{A}$, on
d\'{e}finit le groupe de Grothendieck scind\'{e}
$\langle\mathcal{A}\rangle$. C'est le groupe libre sur les objets
de $\mathcal{A}$ modulo les relations $M=M'+M''$ chaque fois que
$M\cong M'\oplus M''$. Chaque objet $A\in \mathcal{A}$ d\'{e}finit
un \'{e}l\'{e}ment $\langle A \rangle \, \in \, \langle \mathcal{A}
\rangle$.
\end{defn}

L'alg\`ebre $R$ est gradu\'e de la mani\`ere  suivante : $R=\bigoplus_{i\in
\mathbb{Z}}R_i$ avec $R_2 = V^*$ et $R_i=0$ pour $i$ impair. Soit
$\Re=\Re_V$ la cat\'{e}gorie des $R$-bimodules
$\mathbb{Z}$-gradu\'{e}s qui sont de type fini \`{a} gauche et
\`{a} droite. Le groupe $\langle \Re \rangle$ est un anneau
pour $\otimes_R.$

\begin{defn}
Pour chaque objet gradu\'{e} $M=\bigoplus_i M_i,$ et chaque entier $n$, on
 d\'{e}finit l'objet
d\'{e}cal\'{e} $M(n)$ par $(M(n))_i=M_{i+n}.$
\end{defn}

On note $R^s$ le sous-anneau de $R$ des invariants pour
l'action de $s\in W$. Ces d\'efinitions permettent de
formuler le th\'{e}or\`{e}me fondamental de Soergel (\cite{S1}, thm.
1.10) :

\begin{thm}[Soergel]
Il existe un et seulement un morphisme d'anneaux $\varepsilon :
\mathcal{H}\rightarrow \langle \Re\rangle$ tel que $\varepsilon(v)
=\langle R(1) \rangle$ et $\varepsilon(T_s +1) =\langle R\otimes_{R^s}R
\rangle$, pour tout $s\in \mathcal{S}.$
\end{thm}

\begin{defn}
La cat\'{e}gorie $\mathbf{B}$ des bimodules de Soergel est la sous
cat\'{e}gorie de $\Re$ dont les objets sont les $B\in \Re$ avec
$\langle B \rangle$ dans l'image de $\varepsilon$.
\end{defn}

C'est cette cat\'{e}gorie \textbf{B} qu'on \'etudie dans
l'article.

\subsection{Conjecture de positivit\'{e} de Kazhdan-Lusztig} La
 conjecture de positivit\'{e} de Kazhdan-Lusztig pr\'evoit que les coefficients
 des polynômes de Kazhdan-Lusztig sont positifs.  Kazhdan et Lusztig
 l'ont d\'{e}montr\'{e}e pour  $W$ un groupe de Weyl fini ou affine dans
\cite{KL}. Haddad \cite{H} a montr\'{e} cette conjecture pour d'autres
 groupes de
Coxeter finis et Dyer \cite{D} pour les groupes de Coxeter universels.

\subsection{Conjecture de Soergel}

\begin{conj}[Soergel]\label{cs}
Pour tout $x\in W$, il existe un $R$-bimodule ind\'{e}composable
$\mathbb{Z}$-gradu\'{e} $B_x\in \Re$ tel que $\varepsilon(C'_x)=
<B_x>$, o\`u $C'_x$ est l'\'el\'ement de la base de Kazhdan-Lusztig associ\'e \`a $x$.
\end{conj}

\begin{rem}
Dans \cite{S3}, Soergel montre que cette conjecture implique la
conjecture de positivit\'{e} des polyn\^{o}mes de Kazhdan-Lusztig en
construisant un inverse \`{a} gauche de $\varepsilon$.  Cette
 construction est rappel\'ee ci-dessous.
\end{rem}

\begin{rem}
 Dans le cas o\`{u}
$k=\mathbb{C}$ et $W$ est un groupe de Weyl fini, la conjecture de Soergel
 est
montr\'{e}e dans l'article \cite{S1}. Dans \cite{F} Fiebig montre cette
 conjecture pour $W$ un groupe de Coxeter Universel. Dans \cite{S2}
 Soergel montre que si la
caract\'{e}ristique de $k$ est plus grande que le nombre de Coxeter de $W$ et si $W$ est un groupe de Weyl
fini, alors la conjecture de Soergel est \'{e}quivalente \`{a} une
 partie
d'une conjecture de Lusztig portant sur les caract\`{e}res des
repr\'{e}sentations irr\'{e}ductibles de groupes alg\'{e}briques sur
$k$ (par exemple $\mathrm{GL}_n(\mathbb{\bar{F}}_p)$).
\end{rem}

\section{Homomorphismes dans $\mathbf{B}$}
\subsection{}

Nous avons d\'{e}j\`{a} expos\'{e} nos motivations pour \'{e}tudier
la cat\'{e}gorie $\mathbf{B}$. Maintenant posons le probl\`{e}me qui
nous int\'eresse.

Pour chaque $s\in \mathcal{S}$ on note $\theta_s =R\otimes_{R^s}R$.
On va fixer par la suite une \'{e}quation $x_s\in V^*$ de
l'hyperplan des points fixes par $s$ ($x_s$ est bien d\'{e}fini \`{a}
scalaire pr\`{e}s). Dans cet article $\mathrm{Hom}$ va \^{e}tre
 l'espace de
morphismes de $R$-bimodules :
$$\mathrm{Hom}(M,N)=\mathrm{Hom}_{R\otimes R}(M,N)$$

On a que $\mathrm{Hom}(M,N)$ est un $(R,R)-$bimodule $\mathbb{Z}$-gradu\'e. L'action de $R$ \`{a} gauche et \`{a} droite de $\mathrm{Hom}(M,N)$ vient de l'action \`{a} gauche et \`{a} droite
 sur $M$,
ou de mani\`{e}re \'{e}quivalente, sur N. En formules,
$(rf)(m)=f(rm)=r(f(m))$, $(fr)(m)=f(mr)=f(m)r  $, pour tout $r\in
R, m\in M, f\in \mathrm{Hom}(M,N).$
 Si les bimodules $M,N$ sont gradu\'{e}s, alors $\mathrm{Hom}$ va avoir
 en
 plus une structure de $R$-bimodule gradu\'{e}, avec
 $$ \mathrm{Hom}(M(\lambda),N(\lambda'))=\mathrm{Hom}(M,N)(\lambda'
 -\lambda)$$
\begin{prob}\label{pr}
D\'ecrire explicitement $\mathrm{Hom}(\theta_{s_1} \cdots\theta_{s_n},
\theta_{t_1}\cdots \theta_{t_k})$ pour   \linebreak $s_1,\ldots , s_n, t_1, \ldots, t_k \in \mathcal{S}.$
\end{prob}
 Soit \textbf{C} la sous-cat\'egorie pleine de \textbf{B}, dont les
 objets sont des sommes directes finies
d'objets de la forme $\theta_{s_1}...\theta_{s_n}(d)$, pour un $d\in
\mathbb{Z}$. Alors si on r\'esoud le probl\`eme \ref{pr}, on trouve les
 espaces de morphismes dans la cat\'egorie \textbf{C}.
 Le lemme suivant montre que $\mathbf{B}$ est l'enveloppe
 Karoubienne de \textbf{C}.

\begin{lem}
Un bimodule gradu\'{e} $B\in \Re$ appartient \`{a} $\mathbf{B}$ si
et seulement s'il existe $C,D \in \mathbf{C}$ tels que
$$ B\oplus C \cong D$$
\end{lem}

\begin{proof} Soit $\overline{s} =(s_1,...,s_m)$ une suite finie
 quelconque de
r\'{e}flexions simples, $b(\overline{s})=(T_{s_1}+1)...(T_{s_m}+1)$
un \'{e}l\'{e}ment dans l'alg\`{e}bre de Hecke, et on d\'{e}finit le
bimodule
$$\theta_{\overline{s}}=R\otimes_{R^{s_1}}R\otimes _{R^{s_2}} R \otimes
 \cdots  R\otimes_{R^{s_m}}R $$

Comme $<\theta_{\overline{s}}[n]>=
\varepsilon(v^nb(\overline{s}))$, alors
$\theta_{\overline{s}}[n]\in \mathbf{B}.$ Ceci montre que notre
crit\`{e}re est suffisant. Comme les $v^nb(\overline{s})$
engendrent $\mathcal{H}$ comme groupe ab\'elien, il est n\'{e}cessaire. \end{proof}

Avec le lemme suivant, le probl\`{e}me \ref{pr} se r\'{e}duit \`{a}
trouver $\mathrm{Hom}(\theta_{s_1}\cdots \theta_{s_n}, R).$

\begin{lem}\label{bijection}
Pour $M,N \in \mathbf{B}$, le morphisme  

\begin{displaymath}
\begin{array}{lll}
\mathfrak{F} :
 \mathrm{Hom}(\theta_s
M,N) &\rightarrow &  \mathrm{Hom}(M, \theta_s N)(2) \\
\ \ \ \ \ \ \ \ \ \ \ \ f &\mapsto & (m
\mapsto x_s\otimes f(1\otimes m)+ 1\otimes f(1\otimes x_s m))     
\end{array}
\end{displaymath}
 est un isomorphisme de $R$-modules \`a droite gradu\'es.
\end{lem}

\begin{proof}
Si $g\in \mathrm{Hom}(M, \theta_s N)(2)$, on peut \'{e}crire de
 mani\`{e}re
unique $g(m)=1 \otimes g_1(m) +x_s \otimes g_2(m)$, avec $g_1(m),
g_2(m)\in N$. Ceci d\'{e}finit les morphismes $g_1$ et $g_2$
associ\'{e}s \`{a} $g$. Soit $\mathfrak{G} : \mathrm{Hom}(M, \theta_s
 N)(2)
\to \mathrm{Hom}(\theta_s M,N)$ le morphisme qui envoie $g$ vers le morphisme
$\lambda \otimes m \mapsto \lambda g_2(m)$, avec $\lambda \in R$ et
 $m\in M$. Il suffit de montrer que
$\mathfrak{F}$ et $\mathfrak{G}$ sont bien d\'{e}finis et inverses
l'un de l'autre.

 On a $g(x_sm)=x_sg(m)$, c'est \`{a} dire,
 $$1\otimes g_1(x_sm) + x_s \otimes g_2(x_sm)= x_s \otimes g_1(m) +
 1\otimes (x_s)^2 g_2(m). $$
 Par unicit\'{e} de la d\'{e}composition, cette formule permet de
 conclure que $g_2(x_sm)=g_1(m)$. Avec cette formule on montre
 directement que $\mathfrak{F}$ et $\mathfrak{G}$ sont inverses l'un de
 l'autre.

 En outre, on a  $g(r^sm)=r^sg(m)$ pour $r^s\in R^s$. Ceci
 revient \`{a}

\begin{displaymath}
\begin{array}{lll}
 1\otimes g_1(r^sm) +x_s\otimes g_2(r^sm)&=& r^s\otimes g_1(m) +r^sx_s
 \otimes g_2(m)  \\

&=& 1\otimes r^s g_1(m) +x_s \otimes r^s g_2(m).
\end{array}
\end{displaymath}

Ceci implique que $g_2(r^sm)=r^sg_2(m)$, et avec ceci on peut voir
que $\mathfrak{G}(g)$ est bien d\'{e}fini. C'est direct de voir que
c'est un morphisme de bimodules. Voir que
$\mathfrak{F}(f)(m\cdot r)=(\mathfrak{F}(f)(m))\cdot r$ pour $r\in R$ est
aussi direct. On veut alors
$\mathfrak{F}(f)(r\cdot m)=r\cdot (\mathfrak{F}(f)(m))$. Pour  $r\in R^s$ c'est
facile, donc, avec notre d\'{e}composition $R=R^s \oplus x_s R^s$ il
reste \`{a} le montrer pour $r=x_s$. Mais ceci d\'{e}coule
directement du fait que $(x_s)^2 \in R^s$. \end{proof}

\subsection{}

\textbf{Notation.} \'Etant donn\'{e} un espace vectoriel
$\mathbb{Z}$-gradu\'{e} $V=\bigoplus_i V_i$, avec dim$(V)<\infty$, on d\'{e}finit sa dimension
gradu\'{e} par
$$ \underline{\mathrm{dim}}V =\sum (\mathrm{dim} V_i)v^{-i}\in \mathbb{Z}[v,v^{-1}].$$

Soit $R^+$ l'id\'eal de $R$ engendr\'e par les \'el\'ements homog\`enes de degr\'e diff\'erent de z\'ero. On d\'efinit le rang gradu\'e d'un  $R$-module \`a droite $\mathbb{Z}$-gradu\'{e} de type fini $M$  par

$$\underline{\mathrm{rk}}M = \underline{\mathrm{dim}}(M/MR^+)\in \mathbb{Z}[v,v^{-1}]. $$

On a alors  $\underline{\mathrm{dim}}(V(1))=v(\underline{\mathrm{dim}}V)$ et
$\underline{\mathrm{rk}}(M(1))=v(\underline{\mathrm{rk}}M).$ On d\'{e}finit
$\underline{\overline{\mathrm{rk}}}M$ comme l'image de $\underline{\mathrm{rk}}M$ par
 $v\mapsto v^{-1}$.

Pour $x\in W$, on d\'{e}finit le $(R,R)-$bimodule $R_x$ comme $R$ avec
 l'action \`{a} droite
tordue par $x$  (en formules : $r\cdot r'=r x(r')$ pour $r\in R_x$ et $r'\in R$), et
 l'action habituelle de $R$ \`{a} gauche.

On rappelle deux r\'{e}sultats de Soergel. Le premier est
le th\'{e}or\`{e}me 5.3 de \cite{S3}, et le deuxi\`{e}me est une partie
 du
th\'{e}or\`{e}me 5.15 du m\^eme article :

\begin{thm}[Soergel]\label{inversa}
Le morphisme $\varepsilon : \mathcal{H}\rightarrow \langle
\Re\rangle$ admet un inverse \`{a} gauche $ \eta : \langle
\Re\rangle \rightarrow \mathcal{H}$ donn\'{e} par
$$ \langle B\rangle \to \sum_{x\in W}
 \underline{\overline{\mathrm{rk}}}\mathrm{Hom}(B, R_x)T_x.$$
\end{thm}

\begin{prop}[Soergel]\label{l}
Si $M, N \in \mathbf{B}$, alors $\mathrm{Hom}(M,N)$ est libre comme
$R$-module \`{a} droite gradu\'{e}.
\end{prop}

\section{Construction de la base des feuilles l\'eg\`eres}
\subsection{}  Pour $x\in
W$, on d\'{e}finit $\mathcal{R}(x)$ comme l'ensemble de toutes les
expressions r\'{e}duites de $x$ (si $l(x)=n$ cet ensemble est un
sous-ensemble de l'ensemble des n-uplets d'\'{e}l\'{e}ments de $\mathcal{S}$). Si
$1$ est l'\'{e}l\'{e}ment unit\'{e} de $W$, on pose $\theta_1 =R$ et
si $\overline{t}=(t_1, \ldots, t_k)\in \mathcal{R}(x)$, on pose
$\theta_{\overline{t}}=\theta_{t_1}\cdots\theta_{t_k}$.

Nous donnons une d\'{e}finition

\begin{defn}
Soit $\tau :\mathcal{H}\rightarrow \mathbb{Z} [v,v^{-1}]$  l'application
d\'{e}finie par
$$\tau\left(\sum_{x\in W}p_xT_x\right) =p_1 \,\,\,\,\,\,\,\,\,\,\,\, (p_x\in \mathbb{Z} [v,v^{-1}]).$$
\end{defn}

Maintenant on peut \'enoncer un corollaire du th\'eor\`eme
 \ref{inversa}.

\begin{cor}\label{grados} Soit $(s_1,\ldots , s_n)\in \mathcal{S}^n.$
On d\'{e}finit les entiers $n_i$ par $\tau((1+T_{s_1})\cdots
(1+T_{s_n})) =\sum_in_iq^i$. Alors, il existe un isomorphisme de
$R-$modules \`{a} droite gradu\'{e}s
$$\mathrm{Hom}(\theta_{s_1}\cdots \theta_{s_n} ,R) \cong \bigoplus_i n_i R(2i). $$

\end{cor}

\begin{proof}
Par la proposition \ref{l}, il existe des entiers $n'_i$ tels que  $\mathrm{Hom}
(\theta_{\overline{s_n}}\ ,R)\cong (\bigoplus_i n'_i R(2i))$ comme $R$-modules \`a droite, et par le
th\'{e}or\`{e}me \ref{inversa}, on a les \'{e}quations

\begin{displaymath}
\begin{array}{lll}
 \tau((1+T_{s_1})\cdots(1+T_{s_n}))&=& \tau \circ \eta \circ
 \varepsilon ((1+T_{s_1})\cdots(1+T_{s_n}))  \\
  &=&\tau \circ \eta (\langle \theta_{s_1} \cdots \theta_{s_n}\rangle)\\
  &=&\underline{\overline{\mathrm{rk}}}\mathrm{Hom} (\theta_{s_1}\cdots \theta_{s_n},R)
 \\
 &=& \underline{\overline{\mathrm{rk}}} (\bigoplus_i n'_i R(2i))\\
&=& \sum n'_i v^{-2i} \underline{\overline{\mathrm{rk}}}R\\
&=& \sum n'_i q^i.
\end{array}
\end{displaymath}

Ceci permet de conclure que $n'_i=n_i$. \end{proof}

\begin{prop}\label{gradozero}
Soient $s,r\in \mathcal{S}$, $s\neq r$ avec $m(s,r)<\infty$. La composante de degr\'{e} z\'ero de
$$\mathrm{Hom}(\underbrace{\theta_s\theta_r\theta_s\cdots}_{m(s,r) \, \mathrm{termes}},
\underbrace{\theta_r\theta_s\theta_r\cdots}_{m(s,r) \, \mathrm{termes}})$$ est
un espace vectoriel de dimension 1.
\end{prop}

\begin{proof} Avant de commencer la d\'emonstration,  donnons d'abord une d\'efinition  :

\begin{defn}
On pose $G:=\{ q^m+\sum_{i<m}a_iq^i$ pour certains $a_i\in \mathbb{Z}$\} o\`u
 $m=m(s,r)$.
\end{defn}

Par le lemme \ref{bijection} et le corollaire \ref{grados}, montrer le
 corollaire \ref{gradozero} est \'{e}quivalent \`{a} montrer que

\begin{equation}\label{gm}
\tau(\underbrace{(1+T_s)(1+T_r)(1+T_s)\cdots (1+T_r)}_{2m})\in G_m  .
\end{equation}

 Pour  un entier $k>0$, d\'efinissons $Z_{2k}=
 \underbrace{T_rT_sT_r\cdots}_{k \  \mathrm{termes}}$ , $Z_{2k-1}= \underbrace{T_sT_rT_s\cdots}_{k\ \mathrm{termes}}$ et
soit $Z_0=1$. Dans le lemme suivant deg($p$) est le degr\'e du polyn\^ome
 $p$, et $[-]$ est la fonction partie enti\`ere.

 \begin{lem}\label{tau}
Dans l'alg\`ebre de Hecke on a  pour tout entier $n$ l'\'egalit\'e
 suivante :

$$
 \underbrace{(1+T_s)(1+T_r)(1+T_s)\cdots}_{2n\ \mathrm{termes}}=\sum_{j=0}^{4n-1}p_{j,n}Z_j$$

avec $\mathrm{deg}(p_{j,n}$)$<n-[j/4]$ et $p_{4n-1,n}=1$.
\end{lem}

\begin{proof} Prouvons-le par r\'ecurrence sur $n$. Pour $n=1$ est
 clair. Supposons le lemme vrai pour $n=k$. Les trois \'equations suivantes
 d\'ecoulent des d\'efinitions :

\begin{displaymath}
\begin{array}{lll}
T_sZ_{2k-1} &=&  (q-1)Z_{2k-1} +qZ_{2k-2} \\ 
 T_rZ_{2k}&=& (q-1)Z_{2k} +qZ_{2k-3} \\ 
T_sT_rZ_{2k} &=& (q-1)Z_{2k+1}+ q(q-1)Z_{2k-3} + q^2Z_{2k-4}.
\end{array}
\end{displaymath}

Ces \'equations peuvent se r\'e\'ecrire de la mani\`ere suivante :

\begin{displaymath}
\begin{array}{llll}
 T_sZ_{j}&=& (q-1)Z_{j} +qZ_{j-1} & \mathrm{si} \,\,\,
 j\,\,\, \mathrm{est}\,\,\, \mathrm{impair}  \\ 
T_rZ_{j}&=& (q-1)Z_{j} +qZ_{j-3} &  \mathrm{si} \,\,\,j
 \,\,\,\mathrm{est}\,\,\, \mathrm{pair}  \\ 
 T_sT_rZ_{j}&=&  (q-1)Z_{j+1}+ q(q-1)Z_{j-3} + q^2Z_{j-4}  & \mathrm{si}\,\,\,j \,\,\,\mathrm{est}\,\,\, \mathrm{pair}  .
\end{array}
\end{displaymath}

Si on utilise ces \'equations dans le d\'eveloppement du terme \`a
 droite de l'\'egalit\'e suivante
\begin{equation}
 \underbrace{(1+T_s)(1+T_r)(1+T_s)\cdots}_{2(k+1)\ \mathrm{termes}}=(1+T_s+T_r+T_sT_r)\sum_{j=0}^{4k-1}p_{j,k}Z_j
\end{equation}
 on arrive pour tout $j$ \`a exprimer   explicitement $p_{j,k+1}$ en
 fonction de l'ensemble $\{p_{j,k}\}_j$ et de $q$. Les in\'egalit\'es de
 l'\'enonc\'e d\'ecoulent des m\^emes in\'egalit\'es pour les $p_{j,n}$,
 de mani\`ere routini\`ere, et la deuxi\`eme assertion d\'ecoule du fait
 que $p_{4(k+1)-1,k+1}=p_{4k-1,k}$. \end{proof}

\textbf{Preuve de la prop. \ref{gradozero}}. On a (cf \cite[proposition 8.1.1]{GP})

\begin{equation} \tau(T_xT_{y^{-1}})=  \begin{cases} q^{l(x)}
 \text{ si } x=y \\
 0 \,\,\,\,\,\,\,\,\,\, \text{si} \,\,\, x \neq y
 \end{cases}\end{equation}

Ceci implique que

\begin{equation}
\tau(Z_j)= \begin{cases} 0 \,\,\,\, \text{ si } \,\, 0<j<4m-1 \\   q^m
 \text{ si } \,\, j=4m-1
\end{cases}
\end{equation}

Donc si on applique $\tau$ aux deux c\^ot\'es de l'\'egalit\'e du lemme
 \ref{tau}, \'etant donn\'e par ce lemme que deg($p_{0,n}$)$<n$ on
 obtient bien (\ref{gm}). \end{proof}

 Pour chaque couple $(s,r)$ comme dans la proposition \ref{gradozero}, on
choisit  un \'{e}l\'{e}ment non nul $f_{s,r}$ de $$\mathrm{Hom}(\underbrace{\theta_s\theta_r\theta_s\cdots}_{m(s,r) \, \mathrm{termes}},
\underbrace{\theta_r\theta_s\theta_r\cdots}_{m(s,r) \, \mathrm{termes}}).$$  Par le corollaire pr\'{e}c\'{e}dent, il est
bien d\'{e}fini \`{a} scalaire pr\`{e}s. Dans la sous-section \ref{pio},
on va \'{e}liminer cette incertitude, c'est \`{a} dire pour chaque
$s$ et $r$ on va avoir un morphisme $f_{s,r}$ bien d\'{e}fini.

\subsection{}
On sait que pour $s\in \mathcal{S}$ on a  $R \cong R^s \oplus x_s R^s$
comme $R^s$-module \`a gauche. Ceci nous dit que $\theta_s$ admet
$\{1\otimes 1, 1\otimes x_s\}$ comme base en tant que R-module
\`a gauche. Et plus encore, si, pour $s\in \mathcal{S}$, on d\'{e}finit
$x_s^0=1$ et $x_s^1=x_s$, alors  par r\'ecurrence on voit que si
$(t_1,\ldots, t_r)$ est un $r$-uplet d'\'{e}l\'{e}ments de
$\mathcal{S}$, alors $\theta_{t_1} \cdots \theta_{t_r}$ admet $$\{
1\otimes x_{t_1}^{i_1} \otimes \cdots \otimes x_{t_r}^{i_r}
\}_{(i_1,\ldots , i_r) \in \{0,1\}^r}$$ comme base en tant que
R-module \`a gauche.
\begin{defn}\label{pang}
On appelle cette base la \og base normale \fg \  de $\theta_{t_1} \cdots
\theta_{t_r}$. On d\'{e}finit $1\otimes x_{t_1} \otimes x_{t_2}
\otimes \cdots \otimes x_{s_r}$ (resp. 1) comme \og l'\'{e}l\'{e}ment normal \fg \  de
$\theta_{t_1}\cdots \theta_{t_r}$ (resp. $R$). Si $x\in \theta_{t_1} \cdots
\theta_{t_r}$, on va appeler \og partie normale de $x$\fg \  le coefficient
de l'\'{e}l\'{e}ment normal dans la d\'{e}composition de $x$ dans la
base normale.

\end{defn}

\subsection{}\label{pio}
\begin{lem}
Soient $s\neq r\in \mathcal{S}$ avec $m(s,r)\neq \infty.$ On pose $m=m(s,r)$ et
 $$X:=\underbrace{\theta_s\theta_r\theta_s\cdots}_{m \ \mathrm{termes}} \,\ \  et \  \ \, X':=\underbrace{\theta_r\theta_s\theta_r\cdots}_{m\  \mathrm{termes}}$$ Soit $x$ l'\'el\'ement normal de $X$.  Il existe un unique morphisme $f_{s,r}\in \mathrm{Hom}(X,X')$ telle que la partie normale de $f_{s,r}(x)$ soit \'egale \`a $1.$
\end{lem}

\begin{proof}
 Dans \cite{S3}, Soergel montre que si les
$B_x$ existent (voir conjecture 1.14), ils sont uniques \`{a}
isomorphisme pr\`{e}s, et il montre cette conjecture pour un
groupe di\'edral fini. Soit $W'$ le sous-groupe de $W$ engendr\'e par $r$ et $s$. Soit $w_0$ le plus long \'{e}l\'{e}ment de $W'$ dans l'ordre de Bruhat. Dans \cite{S1} Soergel montre que $B_{w_0}=R\otimes_{R^{W'}}R$. Par \cite[ch. IV, cor. 1.11 a.]{Hi} on a  :
$$R \cong \bigoplus_{w\in W} R^{W'}(-2l(w))  \textrm{ comme } 
 R^{W'}-\mathrm{mod} \textrm{ gradu\'e \`a gauche}  $$
et ceci implique

$$B_{w_0}\cong \bigoplus_{w\in W} R(-2l(w)) \textrm{ comme } 
 R-\mathrm{mod} \textrm{ gradu\'e \`a gauche}. $$

Si $M$ est un R-module, on  note $\overline{M}=M/R^+M$. La
derni\`{e}re ligne montre que

$$   \overline{B_{w_0}}\cong \bigoplus_{w\in W} k(-2l(w))  
\textrm{ comme }  k-\mathrm{ev}  \textrm{ gradu\'e}. $$

En particulier,  $\overline{B_{w_0}}$ est de
dimension 1 comme espace vectoriel en degr\'{e} $2m$ :

\begin{equation}\label{1}
(\overline{B_{w_0}})_{2m} \cong k.
\end{equation}

On voit en outre que  $$\{
\overline{1\otimes x_{s}^{i_1} \otimes  x_{r}^{i_2} \otimes
x_s^{i_3} \otimes \cdots} \}_{(i_1,\ldots , i_m) \in \{0,1\}^m}$$
est une base de $\overline{X}$ comme $k$-espace vectoriel gradu\'e (voir d\'efinition \ref{pang}). En
particulier on a

\begin{equation}\label{2}
(\overline{X})_{2m} \cong k
\end{equation}

parce que l'\'{e}l\'{e}ment normal de X est le seul de la base
normale dont le degr\'{e} est $2m$. Dans \cite[prop. 6.16]{S3}, Soergel montre l'existence d'isomorphismes de $(R,R)$-bimodules gradu\'es :

\begin{equation}\label{3}
\mu : X \xrightarrow{\sim} B_{w_0}\oplus M
\end{equation}
et
\begin{equation}\label{4}
\nu : X' \xrightarrow{\sim} B_{w_0}\oplus M'
\end{equation}
pour certains bimodules $M$, $M'$. En passant au quotient on obtient un isomorphisme de k-espaces vectoriels gradu\'es

\begin{equation}\label{5}
\overline{\mu} :    \overline{X}\xrightarrow{\sim} \overline{B_{w_0}}\oplus
    \overline{M}.
\end{equation}

Rappelons que $x$ est l'\'{e}l\'{e}ment normal de $X$. Soit sa
d\'{e}composition $\mu(x)=x_1+x_2$ comme dans (\ref{3}). Par (\ref{1}), (\ref{2})
 et
(\ref{5}), on voit que $\overline{\mu}(\overline{x})=\overline{\mu(x)} \in \overline{B_{w_0}}$. Ceci dit
que  $x_2\in R^+M \subset R^+(B_{w_0}\oplus M)$. Comme $x\notin R^+X$, alors $\mu(x)\notin R^+(B_{w_0}\oplus M),$ parce que $\mu$ est un isomorphisme de $R$-modules. Donc on a

\begin{equation}\label{papa}
x_1=\mu(x)-x_2\notin R^+(B_{w_0}\oplus M).
\end{equation}

En utilisant les identifications (\ref{3}) et (\ref{4}), on d\'efinit $\mathcal{K} \in \mathrm{Hom}(X,X')$,  identit\'{e} sur $B_{w_0}$ et z\'{e}ro sur $M$. Par la proposition \ref{gradozero},  $\mathcal{K}$ est l'unique morphisme de degr\'{e}
z\'{e}ro,  \`{a} scalaire pr\`{e}s, de $\mathrm{Hom}(X,X')$.

On va montrer maintenant que la partie normale de $\mathcal{K}(x)$
est non nulle. Supposons qu'elle soit nulle. Comme $\mathcal{K}$ est un morphisme gradu\'{e} de degr\'{e}
z\'{e}ro, $\mathcal{K}(x)$ est de degr\'{e} $2m$ et comme tous les
\'{e}l\'{e}ments de la base normale sauf l'\'{e}l\'{e}ment normal
sont de degr\'{e} inf\'{e}rieur \`{a} $2m$, on a
$\mathcal{K}(x)\in R^+X'$. Mais d'autre part, $\mathcal{K}(x)=\nu^{-1}(x_1)$. Comme $\nu^{-1}$ est un isomorphisme de $R-$modules, (\ref{papa})  nous dit que $\nu^{-1}(x_1)\notin R^+X'$, ce qui nous donne une
contradiction.

Finalement, comme la partie normale de $\mathcal{K}(x)$ est non
nulle, on choisit $f_{s,r}$ comme le multiple de $\mathcal{K}$ tel
que la partie normale de $f_{s,r}(x)$ soit $1$. \end{proof}

\subsection{}

Pour chaque $s\in \mathcal{S}$ on va d\'{e}finir six morphismes. Consid\'erons la
d\'{e}composition $R=R^s \oplus x_s R^s$. Soient $p_1, p_2\in R^s$ et $p,q,r\in R$. On d\'efinit des morphismes de $(R,R)-$bimodules gradués :

\begin{displaymath}
\begin{array}{lll}\smallskip
 P_s : R \rightarrow R,&&  p_1+x_sp_2 \mapsto p_1  \\ \smallskip
 I_s : R \rightarrow R,&&  p_1+x_sp_2 \mapsto x_sp_2\\ \smallskip
 I'_s : R(2) \rightarrow R,&& p_1+x_sp_2 \mapsto p_2\\ \smallskip
 m^s : \theta_s\rightarrow R, && R\otimes_{R^s}R\ni p\otimes q \mapsto pq\\ \smallskip
 i^s_0 : \theta_s\theta_s (2)\rightarrow R, && R\otimes_{R^s}R\otimes_{R^s}R \ni p\otimes q\otimes r \mapsto pI'_s(q)r\\ \smallskip
 i^s_1 : \theta_s\theta_s (2)\rightarrow \theta_s, && R\otimes_{R^s}R\otimes_{R^s}R \ni p\otimes q\otimes r \mapsto pI'_s(q)\otimes r \in R\otimes_{R^s}R.
\end{array}
\end{displaymath}

Il est facile de voir que ces morphismes sont bien d\'efinis, parce que $I'_s(r^sp)=r^sI'_s(p)$ pour $r^s\in R^s.$

\subsection{}\label{para}
On va fixer jusqu'\`a la fin de la section \textbf{5} une suite
 $(s_1,s_2, \ldots)$ d'\'{e}l\'{e}ments de
$\mathcal{S}$. On d\'{e}finit $\overline{s_0}=1$ et
$\overline{s_n}=(s_1,s_2, \ldots, s_n)$, pour $n\geq 1$.
On sait que si $x\in W$, $s\in \mathcal{S}$ et $l(xs)<l(x)$, alors il
existe une expression r\'{e}duite de $x$ ayant $s$ comme dernier
\'{e}l\'{e}ment. Dans \cite[ch.4, §1, prop. 4]{B}, on montre qu'on
peut passer d'une expression r\'{e}duite  d'un \'{e}l\'{e}ment
\`{a} n'importe quelle autre par une suite de mouvements de
tresse. Donc pour chaque couple $(n,\overline{t})$, avec $n\in
 \mathbb{N}$, $\overline{t}=(t_1,\ldots
 ,t_k)$ une expression r\'eduite de $x:=t_1\cdots t_k \in W$ et avec $l(x s_n) < l(x)$, l'ensemble de suites d'\'el\'ements de
 $\mathcal{R}(x)$
 $$((t^1_1,\ldots, t^1_k),(t^2_1,\ldots, t^2_k), \ldots ,(t^l_1,\ldots,
 t^l_k) )$$
o\`{u} $t^1_i=t_i$ pour tout $1\leq i \leq k$, $t^l_k=s_n$, et o\`u on
passe de $(t^i_1,\ldots, t^i_k)$ vers $(t^{i+1}_1,\ldots
,t^{i+1}_k)$ par un mouvement de tresses, est non vide.

Alors on
 choisit arbitrairement et jusqu'\`a la fin de la section \textbf{5}, pour
 chaque couple $(n,\overline{t})$  un \'el\'ement de cet ensemble, qu'on
 appellera $P(n,\overline{t})$.

 Soit $P(n,\overline{t})= ((t^1_1,\ldots, t^1_k),(t^2_1,\ldots,
 t^2_k), \ldots ,(t^l_1,\ldots, t^l_k) )$.  Pour chaque $1\leq i\leq l-1$ on
 a un morphisme $\pi_i$
associ\'{e} au mouvement de tresses dans
 $\mathrm{Hom}(\theta_{t^i_1} \cdots
\theta_{t^i_k}, \theta_{t^{i+1}_1}\cdots \theta_{t^{i+1}_k})$, du
type $\mathrm{Id}^p\otimes f_{s,r} \otimes \mathrm{Id}^{k-p-m(s,r)}$ pour un certain
$0\leq p \leq k-m(s,r)$. On d\'{e}finit
 $F_n(\overline{t})=\pi_{l-1}\circ \cdots \circ \pi_2 \circ \pi_1\in
\mathrm{Hom}(\theta_{t_1} \cdots \theta_{t_k}, \theta_{t^l_1}\cdots
\theta_{t^l_{k-1}}\theta_{s_n})$.

 \subsection{}
 On va d\'{e}finir par r\'{e}currence sur $n\geq 0$ un sous ensemble
 $A_n$ de
$$\mathfrak{F_n} :=\coprod_{x\in W}\coprod_{\overline{t}\in \mathcal{R}(x)}
 \mathrm{Hom}(\theta_{\overline{s_n}},
\theta_{\overline{t}}) $$ o\`{u} $\coprod$ est
l'union disjointe.

On pose tout d'abord $A_0=\{\mathrm{Id}:R\rightarrow R\}$. On va
construire maintenant $A_n$ \`{a} partir de $A_{n-1}$. On pose
$$A^0_{n-1}=A_{n-1} \cap \left( \coprod_{x\in W}
 \coprod_{\substack{\overline{t}\in \mathcal{R}(x) \\ l(xs_n)>l(x) }}
 \mathrm{Hom}(\theta_{\overline{s_n}},
\theta_{\overline{t}})\right) $$

et $ A^1_{n-1}= A_{n-1}-A^0_{n-1}.$

 On va d\'{e}finir quatre morphismes, 
$f^j_{i,n} :A^j_{n-1} \rightarrow \mathfrak{F_n}$ avec
$i,j\in\{0,1\}$.

Pour les deux premiers morphismes, soient $a\in
\mathrm{Hom}(\theta_{s_1}\cdots \theta_{s_{n-1}},\theta_{t_1}\cdots
\theta_{t_k}) \in A^0_{n-1}$ et $\overline{t}=(t_1,\ldots ,t_k)$. On
 pose

$$f^0_{0,n}(a):\theta_{s_1}\cdots \theta_{s_n}\xrightarrow{a\otimes
 \mathrm{Id}}\theta_{t_1}\cdots
\theta_{t_k}\theta_{s_n}\xrightarrow{\mathrm{Id}^k\otimes
m^{s_n}}\theta_{t_1}\cdots \theta_{t_k}$$

$$f^0_{1,n}(a):\theta_{s_1}\cdots \theta_{s_n}\xrightarrow{a\otimes
 \mathrm{Id}}\theta_{t_1}\cdots
\theta_{t_k}\theta_{s_n}.$$

Pour les deux derniers, soit $a\in
\mathrm{Hom}(\theta_{s_1}\cdots \theta_{s_{n-1}},\theta_{t_1}\cdots
\theta_{t_k}) \in A^1_{n-1}$ et soit $\overline{t}=(t_1,\ldots ,t_k)$. On a $ F_n(\overline{t})\in  \mathrm{Hom}(\theta_{t_1}\cdots
\theta_{t_k},\theta_{t'_1}\cdots\theta_{t'_{k-1}}\theta_{s_n}) $. On
 pose

\begin{multline*}f^1_{0,n}(a):\theta_{s_1}\cdots \theta_{s_n}\xrightarrow{a\otimes
 \mathrm{Id}}\theta_{t_1}\cdots
\theta_{t_k}\theta_{s_n}\xrightarrow{F_n(\overline{t})\otimes
\mathrm{Id}}\theta_{t'_1}\cdots\theta_{t'_{k-1}}\theta_{s_n}\theta_{s_n}
\\ \xrightarrow{\mathrm{Id}^{k-1}\otimes
i^s_0}\theta_{t'_1}\cdots\theta_{t'_{k-1}}
\end{multline*}

\begin{multline*}f^1_{1,n}(a):\theta_{s_1}\cdots \theta_{s_n}\xrightarrow{a\otimes
 \mathrm{Id}}\theta_{t_1}\cdots
\theta_{t_k}\theta_{s_n} \xrightarrow{F_n(\overline{t})\otimes
\mathrm{Id}}\theta_{t'_1}\cdots\theta_{t'_{k-1}}\theta_{s_n}\theta_{s_n}\\ \xrightarrow{\mathrm{Id}^{k-1}\otimes
i^s_1}\theta_{t'_1}\cdots\theta_{t'_{k-1}}\theta_{s_n}.
\end{multline*}

Maintenant on peut d\'{e}finir $A_n$

$$ A_n=\bigcup_{0\leq i,j \leq 1}f^j_{i,n}(A^j_{n-1}). $$

On d\'{e}finit $A'_n$ comme le sous ensemble de $A_n$ des
\'{e}l\'{e}ments appartenant \`{a}
 $\mathrm{Hom}(\theta_{\overline{s_n}},
R).$

\section{Le th\'{e}or\`{e}me et sa preuve}

\begin{thm}\label{bacan}
L'ensemble $A'_n$ est une R-base de
 $\mathrm{Hom}(\theta_{\overline{s_n}}, R).$
\end{thm}

\begin{rem}
On va appeler $A'_n$ une \textsl{\og base de feuilles l\'eg\`eres\fg} : En
 tensorisant par l'identit\'e on peut voir $A_k \subset \mathfrak{F_n}$
 pour $k\leq n$. Alors on peut voir $A_n$ comme un arbre binaire parfait,
 dont on peut associer \`a chaque feuille (une feuille est un morphisme
 de $\theta_{s_1}\cdots \theta_{s_n}$ vers $\theta_{t_1}\cdots
 \theta_{t_k}$)  un nombre ou \og poids\fg \  (dans ce cas le poids serait $k$). Avec ce
 point de vue, $A'_n$ est l'ensemble de feuilles qui ont poids z\'ero. 
\end{rem}

\begin{rem}
Il faut noter que chaque morphisme de $A'_n$ dépend du choix des  $P(n,\overline{t})$. On va donner un exemple  dans lequel différents choix du $P(n,\overline{t})$ donnent des morphismes associés différents. 

Soient $s,r\in \mathcal{S}$ avec $m(s,r)=3$. Soient $$f=(m^s \circ i_1^s)\circ (Id \otimes m^r \circ i_1^r \otimes Id)\circ (Id^2 \otimes m^s \circ i_1^s \otimes Id^2)$$ et $$g=f\circ (f_{r,s}\otimes Id^3)\circ(f_{s,r}\otimes Id^3)$$
  deux morphismes appartenant à $\mathrm{Hom}(\theta_s\theta_r\theta_s\theta_s\theta_r\theta_s,R)$. On va montrer que $f\neq g.$ Pour ceci on définit $$\bar{x}= 1\otimes x_r\otimes 1 \otimes x_s \otimes 1 \otimes x_s \otimes 1 + 1\otimes 1\otimes  x_r\otimes x_s \otimes 1 \otimes x_s \otimes 1 \in \theta_s\theta_r\theta_s\theta_s\theta_r\theta_s.$$  On voit facilement que $f(\bar{x})=1$. Il nous reste a montrer que $g(\bar{x})=0$. Pour ceci il suffit de montrer que 
\begin{equation}\label{cauro}
 f_{s,r}(1\otimes x_r \otimes 1 \otimes 1+1\otimes 1 \otimes x_r\otimes 1)=0.
\end{equation}

Comme $C'_sC'_rC'_s=C'_{srs}+C'_s$, le théorème 4.2 de \cite{S3} et le lemme 2 de \cite{S1} impliquent que  $\theta_s\theta_r\theta_s(3)\simeq   R\otimes_{R^{\langle s,r \rangle}}R(3) \oplus \theta_s(1)$, c'est-à-dire, 
\begin{equation}\label{pacho}
 \theta_s\theta_r\theta_s\simeq   R\otimes_{R^{\langle s,r \rangle}}R \oplus \theta_s(-2).
\end{equation}
Il est facile de voir que, à scalaire près, le seul morphisme de degré zéro de $R\otimes_{R^{\langle s,r \rangle}}R$ vers $\theta_s\theta_r\theta_s$
est le morphisme $i_{s,r}$ défini par $(1\otimes 1\mapsto 1\otimes 1\otimes 1\otimes 1).$ Par le corollaire \ref{grados} on voit qu'à scalaire près il y a un seul morphisme de degré zéro de $\theta_s(-2)$ vers $\theta_s\theta_r\theta_s$. Ce morphisme est le morphisme $(1\otimes 1\mapsto 1\otimes x_r\otimes 1 \otimes 1+1\otimes 1\otimes x_r \otimes 1)$. Donc on a (toujours à scalaire près) explicité l'isomorphisme (\ref{pacho}).
 
Comme on a vu dans le lemme \ref{pio}, $R\otimes_{R^{\langle s,r \rangle}}R\simeq B_{srs}$, donc on a un diagramme commutatif
$$\xymatrix{
{\theta_s\theta_r\theta_s}\ar[r] \ar[dr]_{f_{s,r}} &{R\otimes_{R^{\langle s,r \rangle}}R}\ar[d]^{i_{r,s}}\\
{}&{\theta_r\theta_s\theta_r}
}
$$
 la flèche horizontale étant la surjection canonique en (\ref{pacho}). Ceci permet de prouver (\ref{cauro}), et donc de conclure que $f\neq g.$
\end{rem}

\textbf{Exemple :} Dans l'exemple suivant nous avons
 $(s_1,s_2,s_3,s_4)=(s,r,s,r)$ et $m(s,r)=3$. Par des raisons d'espace nous montrons
 seulement une moiti\'e de l'arbre, la \og moiti\'e gauche\fg, c'est \`a dire, nous ne
 montrons pas les morphismes qui passent par l'\'el\'ement encadr\'e  $\theta_r
 \theta_s \theta_r$. Nous avons encercl\'e en bas \`a droite la seule
 feuille l\'eg\`ere de cette moiti\'e d'arbre.

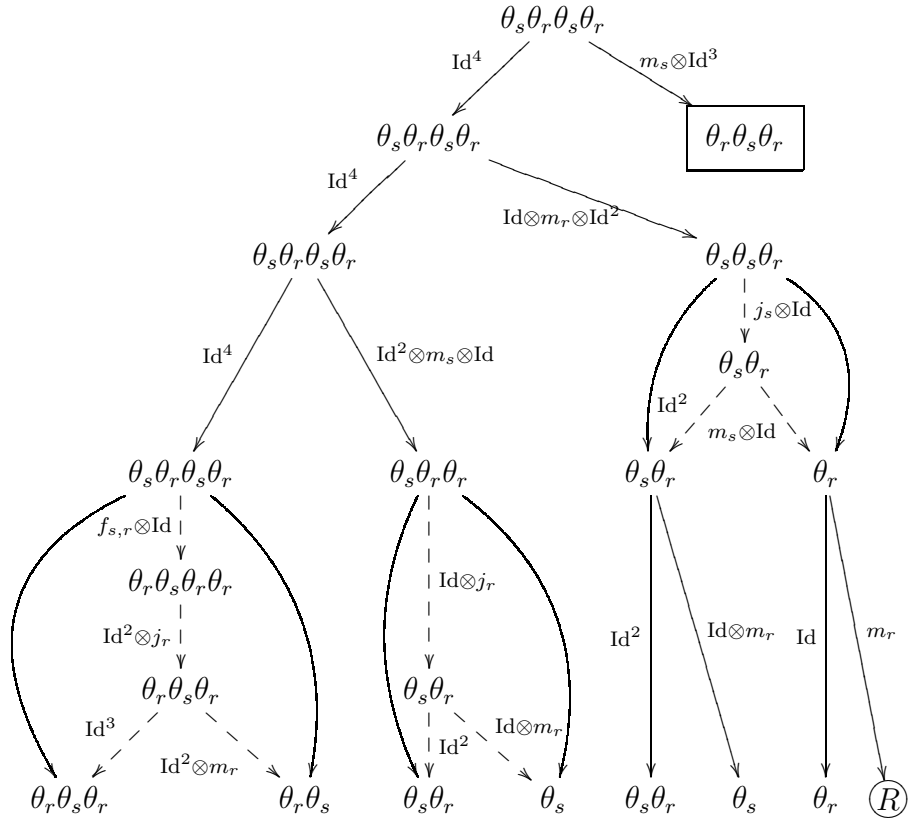
\xymatrix@C=0cm{
&&&&\theta_s \theta_r \theta_s\theta_r \ar[dl]_{\mathrm{Id}^4} \ar[drr]^{m_s
 \otimes \mathrm{Id}^3} & &\\
 & &&\theta_s \theta_r \theta_s\theta_r \ar[dl]_{\mathrm{Id}^4}
 \ar[drrr]_{\mathrm{Id}\otimes m_r \otimes \mathrm{Id}^2} & & & *+<0.5cm>[F]{\theta_r \theta_s\theta_r} &\\
 &&\theta_s \theta_r \theta_s\theta_r \ar[ddl]_{\mathrm{Id}^4}
 \ar[ddr]^{\mathrm{Id}^2\otimes m_s \otimes \mathrm{Id}} & &&&\theta_s \theta_s\theta_r
 \ar@{-->}[d]^{j_s\otimes \mathrm{Id}}\ar@/_5mm/[ddl] \ar@/^8mm/[ddr]\\
 &&&&& &\theta_s \theta_r \ar@{-->}[dl]_{\mathrm{Id}^2} \ar@{-->}[dr]_{m_s\otimes
 \mathrm{Id}} \\
& \theta_s \theta_r \theta_s\theta_r \ar@/^8mm/[dddr] \ar@/_15mm/[dddl]
  \ar@{-->}[d]_{f_{s,r}\otimes \mathrm{Id}}&&   \theta_s \theta_r \theta_r
  \ar@/^9mm/[dddr] \ar@/_6mm/[ddd]\ar@{-->}[dd]^{\mathrm{Id}\otimes j_r} & &\theta_s
 \theta_r \ar[ddd]_{\mathrm{Id}^2} \ar[dddr]^{\mathrm{Id}\otimes m_r}& & \theta_r \ar[ddd]_{\mathrm{Id}}
 \ar[dddrr]^{m_r}  \\
& \theta_r \theta_s \theta_r\theta_r \ar@{-->}[d]_{\mathrm{Id}^2\otimes j_r}
  \\
&    \theta_r \theta_s \theta_r \ar@{-->}[dl]_{\mathrm{Id}^3}
 \ar@{-->}[dr]_{\mathrm{Id}^2\otimes m_r}&& \theta_s \theta_r   \ar@{-->}[d]^{\mathrm{Id}^2}
 \ar@{-->}[dr]^{\mathrm{Id}\otimes m_r} \\
  \theta_r \theta_s\theta_r && \theta_r \theta_s &\theta_s \theta_r
 &\theta_s &\theta_s \theta_r & \theta_s &  \theta_r & & *+[o][F]{R}  }

\begin{proof}  D\'efinissons les polyn\^omes $p^x_n$ par
 $$(1+T_{s_1})\cdots (1+T_{s_{n}})=
\sum_{x\in W}p^x_nT_x$$

On a les \'{e}quations suivantes
\begin{multline*}
(1+T_{s_1})\cdots (1+T_{s_{n-1}})(1+T_{s_{n}})= \left( \sum p^x_{n-1}T_x\right)+
\\
+\left( \sum_{l(xs_n)>l(x)} p^x_{n-1}T_{xs_n}\right)+\left(\sum_{l(xs_n)<l(x)}
p^x_{n-1}T_xT_{s_n}\right)
\end{multline*}

et

\begin{multline*}
\sum_{l(xs_n)<l(x)} p^x_{n-1}T_xT_{s_n}= \left(\sum_{l(xs_n)<l(x)} q
p^x_{n-1}T_{xs_n}\right)+ \\
+\left(\sum_{l(xs_n)<l(x)}
p^x_{n-1}\underbrace{T_{xs_n}T_{s_n}}_{T_x}(q-1)\right),
\end{multline*}

 qui impliquent

\begin{multline*}
(1+T_{s_1})\cdots (1+T_{s_{n-1}})(1+T_{s_{n}})= \left(\sum_{l(xs_n)>l(x)}p^x_{n-1}T_x +p^x_{n-1}T_{xs_n}\right)+ \\
 +\left(\sum_{l(xs_n)<l(x)} q(p^x_{n-1}T_x+p^x_{n-1}T_{x{s_n}})\right).
\end{multline*}

Cette derni\`{e}re formule, le corollaire \ref{grados}, et la
 construction
r\'ecursive de $A'_n$ vont nous montrer que les degr\'{e}s
gradu\'{e}s des \'{e}l\'{e}ments de $A'_n$ sont ceux que devrait
avoir une base. Avec les d\'{e}finitions suivantes on va dire ceci
d'une mani\`{e}re plus pr\'{e}cise.

\begin{defn}
Soit $X= \{x_1, \ldots,  x_n  \}$ un ensemble d'\'el\'ements homog\`enes
 d'espaces vectoriels gradu\'{e}s. On
d\'{e}finit $\mathbb{Y}(X)=\sum_i
q^{\mathrm{deg}(x_i)/2}.$
\end{defn}

\begin{defn}
Pour $x\in W$ on d\'{e}finit
$$A_n(x)= A_n \cap \coprod_{\overline{t}\in \mathcal{R}(x)}
 \mathrm{Hom}(\theta_{\overline{s_n}},
\theta_{\overline{t}}).$$
\end{defn}

\begin{lem}
On a l'\'egalit\'e
$\mathbb{Y}(A_n(x))=p^x_{n} $. En particulier, $\mathbb{Y}(A'_n)=p^1_n.$

\end{lem}
 \begin{proof} On va le montrer par r\'{e}currence sur $n$.
Pour $n=1$ c'est clair. Supposons-le vrai pour $n-1$. Par construction,
 \`{a} chaque $a\in (A_{n-1}(x))^0$ on associe les
deux \'{e}l\'{e}ments  $a'\in A_n(x)$ et $a''\in A_n(xs_n)$, et
\`{a} chaque $b\in (A_{n-1}(x))^1$ on associe les deux
\'{e}l\'{e}ments  $b'\in A_n(x)$ et $b''\in A_n(xs_n)$, avec
deg($a$)=deg($a'$)=deg($a''$), et deg($b$)$+2$=deg($b'$)=deg($b''$),
ce qui permet de conclure la r\'{e}currence. \end{proof}

Supposons que l'on ait montr\'{e} que les \'{e}l\'{e}ments de $A'_n$
sont lin\'{e}airement ind\'{e}pendants pour  l'action de R. Soit $T$ le
 sous $R$-module  de $\mathrm{Hom}(\theta_{\overline{s_n}}\ ,R)$   engendr\'{e} par les \'{e}l\'{e}ments de $A'_n$.
  Dans chaque degr\'{e}, $T$ et  $\mathrm{Hom}(\theta_{\overline{s_n}}\ ,R)$ sont de
dimension finie comme $k$-espaces vectoriels, et ils ont la m\^{e}me
dimension parce que $\mathbb{Y}(A'_n)=p^1_n$ et $T$ est libre pour
 l'action de R (voir corollaire \ref{grados}), donc ils sont \'{e}gaux. Cela
 \'{e}tant vrai \`{a}
chaque degr\'{e}, on d\'{e}duit que $T=\mathrm{Hom} (\theta_{s_1}\cdots
\theta_{s_n},R)$, et cela finit la preuve du th\'{e}or\`{e}me \ref{bacan}.

 Donc il nous reste 
\`{a} montrer que les \'{e}l\'{e}ments de $A'_n$ sont
lin\'{e}airement ind\'{e}pendants pour l'action de R.

Les \'{e}l\'{e}ments de $A'_n$ sont de la forme $f^{j_n}_{i_n,n}
\circ \cdots \circ f^{j_2}_{i_2,2}\circ f^{j_1}_{i_1,1}(\mathrm{Id})$ avec
$i_k, j_k \in \{0,1\}$. Soit $I_n=\{(i_1,
\ldots i_n)$ tel que $\exists \,\,(j_1, \ldots j_n)$ avec
$f^{j_n}_{i_n,n} \circ \cdots \circ f^{j_1}_{i_1,1}(\mathrm{Id}) \in
A'_n\}\subseteq \{0,1\}^n.$ On remarque que dans cette derni\`{e}re
d\'{e}finition le $n$-uplet $(j_1, \ldots j_n)$ est unique, s'il existe.

On rappelle que $x_s^0=1$ et $x_s^1=x_s $ pour $s\in \mathcal{S}$. Si
$\overline{i}=(i_1,\ldots ,i_n) \in \{0,1\}^n$, on pose
 $x_{\overline{i}}
=1\otimes x^{i_1}_{s_1}\otimes x^{i_2}_{s_2} \otimes \cdots \otimes
x^{i_n}_{s_n} \in \theta_{s_1}\theta_{s_2}\cdots \theta_{s_n}$, et
si $\overline{i}\in I_n$, on pose $f_{\overline{i}}=f^{j_n}_{i_n,n}
\circ \cdots \circ f^{j_1}_{i_1,1}(\mathrm{Id}).$

Sur $\{0,1\}^n$, on appelle $\prec$ l'ordre total suivant : Si $\sum_j
i_j>\sum_j i'_j$ alors $(i_1,\ldots i_n) \succ (i'_1,\ldots i'_n)$.
Sinon, soit $r$ le plus petit entier tel que $i_r \neq i'_r$. Si
$i_r=0$, alors $(i_1,\ldots i_n) \succ (i'_1,\ldots i'_n)$, et si
$i_r=1$, alors $(i_1,\ldots i_n) \prec (i'_1,\ldots i'_n)$.

Pour finir la d\'{e}monstration du th\'{e}or\`{e}me, on veut montrer
que si $\sum_{\overline{i}\in I_n}
a_{\overline{i}}f_{\overline{i}}=0$, avec $a_{\overline{i}} \in R$,
alors $a_{\overline{i}}=0$ pour tout $\overline{i}\in I_n$. Pour
cela il est suffisant de montrer les deux faits suivants :
\begin{enumerate}[(a)]
    \item $f_{\overline{i}}(x_{\overline{i}}) =1 $
    \item $f_{\overline{i}}(x_{\overline{i'}})=0$ pour $\overline{i} \succ
\overline{i'}$
\end{enumerate}

\begin{defn}\label{sup}
On dit qu'un \'{e}l\'{e}ment $x\in \theta_{t_1}\cdots \theta_{t_k}$
est \og sup\'erieur\fg, si $x\in R_+ \theta_{t_1}\cdots \theta_{t_k}$, et il
est \og normalsup \fg s'il appartient \`{a} $1\otimes x_{t_1} \otimes
x_{t_2} \otimes \cdots \otimes x_{s_r}+R_+ \theta_{t_1}\cdots
\theta_{t_k}$. Pour $f\in \theta_{t_1}\cdots \theta_{t_k}$, on note $f\ddag$ l'ensemble $f+R_+
\theta_{t_1}\cdots \theta_{t_k}$.
\end{defn}

 Si $\overline{i}=(i_1,\ldots , i_n)\in I_n$,
on d\'{e}finit $f^m_{\overline{i}} = f^{j_m}_{i_m,m} \circ \cdots
\circ f^{j_2}_{i_2,2}\circ f^{j_1}_{i_1,1}(\mathrm{Id})$, et
$x^m_{\overline{i}}=1\otimes x^{i_1}_{s_1}\otimes x^{i_2}_{s_2}
\otimes \cdots \otimes x^{i_m}_{s_m} \in
\theta_{s_1}\theta_{s_2}\cdots \theta_{s_m}$, pour $1\leq m \leq n$.

\begin{lem}\label{supsup}
Soient $\overline{i}, \overline{i'}\in I_n$. Si
$f^{u-1}_{\overline{i}}(x^{u-1}_{\overline{i'}})$ est sup\'erieur,
 alors
$f^u_{\overline{i}}(x^{u}_{\overline{i'}})$ est sup\'erieur aussi.
\end{lem}

\begin{proof} On a
$f^{u}_{\overline{i}}(x^{u}_{\overline{i'}})=
f^{j_{u}}_{i_{u},u}(f^{u-1}_{\overline{i}})(x^{u-1}_{\overline{i'}}
\otimes x^{i'_{u}}_{s_{u}})$. Par hypoth\`{e}se
$(f^{u-1}_{\overline{i}} \otimes \mathrm{Id}) (x^{u-1}_{\overline{i'}}
\otimes x^{i'_{u}}_{s_{u}})$ est sup\'erieur, et le fait que
$f^u_{\overline{i}}(x^{u}_{\overline{i'}})$ est l'image de cet
\'{e}l\'{e}ment par un morphisme de bimodules, permet de conclure le
lemme. \end{proof}

\begin{lem}\label{nsns}
Le morphisme $F_r(\overline{t})$ envoie toujours un \'{e}l\'{e}ment
normalsup (resp. sup\'erieur) vers un \'{e}l\'{e}ment normalsup (resp.
sup\'erieur).
\end{lem}

\begin{proof}
Comme  $F_r(\overline{t})$ est un morphisme de bimodules, il envoie
un \'{e}l\'{e}ment sup\'erieur vers un \'{e}l\'{e}ment sup\'erieur.

Maintenant on veut montrer que $F_r(\overline{t})$ envoie un
\'{e}l\'{e}ment normalsup vers  un \'{e}l\'{e}ment normalsup. Il
suffit \`{a} son tour de montrer ceci pour les morphismes du type
$\mathrm{Id}^a \otimes f_{s,r} \otimes \mathrm{Id}^b$, ce qui revient \`{a} le montrer
pour les $f_{s,r}$, mais ceci est vrai par d\'{e}finition de
$f_{s,r}.$ \end{proof}

\begin{lem}\label{norsup}
Si $f^m_{\overline{i}} \in  \mathrm{Hom}(\theta_{s_1}\cdots
\theta_{s_m},\theta_{t_1}\cdots \theta_{t_r})$, alors
 $f^m_{\overline{i}}(x^m_{\overline{i}})$ est  un \'{e}l\'{e}ment
normalsup de $\theta_{t_1}\cdots \theta_{t_r}$.
\end{lem}

\begin{proof}: On va le montrer par r\'{e}currence sur $m$. Pour
$m=1$, on sait que $j_1=0$.

Si $i_1=0$, alors $f^0_{0,1}(\mathrm{Id})(1\otimes 1) = 1 \in R$.

Si $i_1=1$, alors $f^0_{1,1}(\mathrm{Id})(1\otimes x_{s_1}) = 1 \otimes
x_{s_1}\in \theta_{s_1}.$

Supposons le lemme vrai pour tout $m<r$. Pour simplifier les
notations on va appeler $g_j:=f^j_{\overline{i}}$. Soit $g_{r-1}\in
\mathrm{Hom}(\theta_{s_1}\cdots \theta_{s_{r-1}},\theta_{t_1}\cdots
\theta_{t_k})$. On note aussi $\overline{t}=(t_1, \ldots t_k)$.
Comme $g_r=f^{j_r}_{i_r,r}(g_{r-1})$, il y a quatre cas :

\begin{enumerate}
  \item $j_r=0$, $i_r=0$, (donc $x^{i_r}_{s_r}= 1 $), alors
\begin{displaymath}
\begin{array}{lll}
g_r(x^r_{\overline{i}}) &=&(\mathrm{Id}^k \otimes m^{s_r})\circ
(g_{r-1}\otimes
  \mathrm{Id})(x^r_{\overline{i}}) \\
&\in & (\mathrm{Id}^r \otimes m^{s_r})(1\otimes x_{t_1} \otimes \cdots \otimes
 x_{t_k}\otimes 1\ddag) \\
&\in& 1\otimes x_{t_1} \otimes \cdots \otimes x_{t_k} \ddag \subseteq
\theta_{t_1} \cdots \theta_{t_k}.
\end{array}
\end{displaymath}

  \item $j_r=0$, $i_r=1$, (donc $x^{i_r}_{s_r}=  x_{s_r}$) alors
  \begin{displaymath}
\begin{array}{lll}
g_r(x^r_{\overline{i}}) &=&(g_{r-1}\otimes
  \mathrm{Id})(x^r_{\overline{i}}) \\
&\in& 1\otimes x_{t_1} \otimes \cdots \otimes x_{t_k}\otimes
x_{t_s}\ddag \subseteq \theta_{t_1} \cdots \theta_{t_k}\theta_{t_s}.
\end{array}
\end{displaymath}

  \item $j_r=1$, $i_r=0$, (donc $x^{i_r}_{s_r}=  1$) alors
  \begin{displaymath}
\begin{array}{lll}
g_r(x^r_{\overline{i}}) &=&(\mathrm{Id}^{k-1}\otimes i^{s_r}_0) \circ
(F_r(\overline{t})  \otimes \mathrm{Id}))\circ (g_{r-1}\otimes
  \mathrm{Id})(x^r_{\overline{i}}) \\
&\in& (\mathrm{Id}^{k-1}\otimes i^{s_r}_0) \circ (F_r(\overline{t})  \otimes
 \mathrm{Id}))(1\otimes x_{t_1} \otimes \cdots \otimes x_{t_k} \otimes 1\ddag) \\
&\in& (\mathrm{Id}^{k-1}\otimes i^{s_r}_0) (1\otimes x_{t'_1} \otimes \cdots
 \otimes x_{t'_{k-1}} \otimes x_{s_r} \otimes 1 \ddag) \\
&\in& 1\otimes x_{t'_1} \otimes \cdots \otimes x_{t'_{k-1}} \ddag\subseteq
\theta_{t'_1} \cdots \theta_{t'_{k-1}}.
\end{array}
\end{displaymath}

  \item $j_r=1$, $i_r=1$, (donc $x^{i_r}_{s_r}=  x_{s_r}$) alors

\begin{displaymath}
\begin{array}{lll}
g_r(x^r_{\overline{i}}) &=&(\mathrm{Id}^{k-1}\otimes i^{s_r}_1) \circ
(F_r(\overline{t})  \otimes \mathrm{Id}))\circ (g_{r-1}\otimes
  \mathrm{Id})(x^r_{\overline{i}}) \\
&\in& (\mathrm{Id}^{k-1}\otimes i^{s_r}_1) \circ (F_r(\overline{t})  \otimes
 \mathrm{Id}))(1\otimes x_{t_1} \otimes \cdots \otimes x_{t_k} \otimes x_{t_r}\ddag)
 \\
&\in& (\mathrm{Id}^{k-1}\otimes i^{s_r}_1) (1\otimes x_{t'_1} \otimes \cdots
\otimes x_{t'_{k-1}} \otimes x_{s_r} \otimes
x_{s_r}\ddag) \\
 &\in& 1\otimes x_{t'_1} \otimes \cdots \otimes x_{t'_{k-1}}\otimes
 x_{s_r}\ddag \subseteq \theta_{t'_1} \cdots
\theta_{t'_{k-1}} \theta_{s_r} .
\end{array}
\end{displaymath}

\end{enumerate}

Le passage de la deuxi\`eme \`a la troisi\`eme ligne dans (3) et (4) vient du fait que
$F_r(\overline{t})$ envoie un \'{e}l\'{e}ment normalsup vers un
\'{e}l\'{e}ment normalsup (lemme \ref{nsns}). \end{proof}

Maintenant on peut finir la preuve du  th\'{e}or\`{e}me. La
premi\`{e}re chose qu'on peut constater est que pour $s \in
\mathcal{S}$, le morphisme $I'_s$ est gradu\'{e} de degr\'{e} $-2$.
Donc $f_{\overline{i}}=f^{j_n}_{i_n,n} \circ \cdots \circ
f^{j_2}_{i_2,2}\circ f^{j_1}_{i_1,1}(\mathrm{Id})$ est un morphisme
gradu\'{e} de degr\'{e} $\sum_{p=1}^n -2j_p.$

Soit $$X_{a,b} = \{p \vert j_p=a\, , \, i_p=b\}. $$

Si $a\in \mathrm{Hom}(\theta_{s_1}\cdots
 \theta_{s_m},\theta_{t_1}\cdots
\theta_{t_k})$, on pose $\varpi(a) = k$. Les quatre
\'{e}quations suivantes sont faciles \`{a} v\'{e}rifier :

\begin{displaymath}
\begin{array}{lll}
\varpi(a) &=&\varpi(f_{0,r}^0(a)) \\
\varpi(a)&=& \varpi(f^1_{1,r}(a)) \\
\varpi(a)+1&=& \varpi(f_{1,r}^0(a)) \\
\varpi(a)-1 &=& \varpi(f_{0,r}^1(a)) .
\end{array}
\end{displaymath}

Si $\overline{i} \in I_n$, alors, comme
$\varpi(f_{\overline{i}})=0$, on peut
d\'{e}duire que $\mathrm{card}(X_{1,0} ) = \mathrm{card}(X_{0,1}) $. Donc on a

\begin{displaymath}
\begin{array}{lll}
\sum_{p=1}^n r_p&=&   \mathrm{card}(X_{1,0})+ \mathrm{card}(X_{1,1}) \\
&=&   \mathrm{card}(X_{0,1})+ \mathrm{card}(X_{1,1})\\
&=& \sum_{p=1}^n i_p.
\end{array}
\end{displaymath}

Donc $f_{\overline{i}}$ est un morphisme gradu\'{e} de degr\'{e}
$\sum_{p=1}^n -2i_p$. On a aussi que si $\overline{i'}=(i'_1, \ldots
, i'_m)$ le degr\'{e} de $x_{\overline{i'}}$ est $\sum 2i'_p$.
 Donc, si $\sum i_p
> \sum i'_p$, alors $f_{\overline{i}}(x_{\overline{i'}})=0$.

Maintenant supposons $\sum i_p = \sum i'_p$. Par le raisonnement
pr\'{e}cedent,
\begin{equation}\label{deg}
\mathrm{deg}(f_{\overline{i}}(x_{\overline{i'}}))= 0,
\end{equation}
c'est \`{a} dire, $f_{\overline{i}}(x_{\overline{i'}})$ est un
scalaire. Ceci et le lemme \ref{norsup} permettent de conclure
la partie a), c'est \`{a} dire, que
$f_{\overline{i}}(x_{\overline{i}}) =1 $.

 Pour montrer la partie b), c'est \`{a} dire que
 $f_{\overline{i}}(x_{\overline{i'}})=0$ pour
$\overline{i} \succ \overline{i'}$, on va supposer $\sum i_p = \sum
i'_p$  , $\overline{i'} \prec \overline{i}$, et aussi
$f_{\overline{i}}(x_{\overline{i'}}) \neq 0$. On va arriver \`{a}
une contradiction.

 Soit $r$ le plus petit entier tel que $i_r \neq
i'_r$, donc $i_r=0$ et $i'_r=1$. Par le lemme pr\'{e}c\'{e}dent,
$f^{r-1}_{\overline{i}}(x^{r-1}_{\overline{i'}})$ est un
\'{e}l\'{e}ment normalsup. Donc si $n$ est l'\'el\'ement normal de l'ensemble d'arriv\'ee du morphisme $f^{r}_{\overline{i}}$,   on a que $f^{r}_{\overline{i}}(x^{r}_{\overline{i'}})$ appartient \`{a}
l'ensemble $(n\cdot x_{s_r})\ddag$ (voir la d\'efinition \ref{sup}).
Le lemme suivant (lemme \ref{ssi}) nous dit que cet \'{e}l\'{e}ment
 est
sup\'erieur, donc le lemme \ref{supsup} et l'\'equation (\ref{deg})
 permettent d'aboutir \`{a} une
contradiction. \end{proof}

\begin{lem}\label{ssi}

Soit $(t_1,\ldots, t_r)$ un r-uplet d'\'{e}l\'{e}ments de
$\mathcal{S}$. Un \'{e}l\'{e}ment $x\in \theta_{t_1}\cdots
\theta_{t_k}$ est sup\'erieur, si et seulement si il est z\'{e}ro, ou
bien s'il existe une \'{e}criture $$x=\sum_{i=1}^m p_0^i\otimes p_1^i
\otimes \cdots \otimes p_k^i$$ qui satisfait la propri\'et\'{e}
suivante, qu'on appellera propri\'et\'e (*M) : les $p^i_j\in R$ sont des
\'{e}l\'{e}ments homog\`{e}nes, et pour tout $1 \leq i \leq m$, il
existe $r_i \in \mathbb{N}_0$, avec $\sum_{j=0}^{r_i} \mathrm{deg}(p_j^i)
\geq r_i +1$. Ici $M=\mathrm{max}\{r_i \vert 1\leq i\leq m\}.$

\end{lem}

\begin{proof} La direction \og seulement si\fg \   est claire parce que si $x$ est sup\'erieur, on a $x=r^+ \cdot z$, avec $r^+\in R^+$ et $z\in \theta_{t_1}\cdots
\theta_{t_k}$. Alors si $z=\sum_{i=1}^m p_0^i\otimes p_1^i
\otimes \cdots \otimes p_k^i$ avec les $p^i_j\in R$  des
\'{e}l\'{e}ments homog\`{e}nes, on a $x=\sum_{i=1}^m r^+p_0^i\otimes p_1^i
\otimes \cdots \otimes p_k^i$, qui satisfait la propri\'et\'e (*0).

Pour l'autre sens, soit $x=\sum_{i=1}^m
p_0^i\otimes p_1^i \otimes \cdots \otimes p_k^i$ satisfaisant la propri\'et\'e (*M). Comme $p^i_{r_i}=
P_{t_{r_i}}(p^i_{r_i})+x_{t_{r_i}}I'_{t_{r_i}}(p^i_{r_i})$, alors

\begin{multline*}
x=\sum_{i=1}^m p_0^i \otimes \cdots \otimes
 p^i_{r_{i}-1}P_{t_{r_i}}(p^i_{r_i})\otimes 1 \otimes \cdots \otimes p_n^i +\\
 +p_0^i \otimes \cdots
 \otimes p^i_{r_{i}-1}I'_{t_{r_i}}(p^i_{r_i})\otimes x_{t_{r_i}}
 \otimes \cdots \otimes p_n^i.
 \end{multline*}

Si $p_0^i \otimes \cdots \otimes
p^i_{r_{i}-1}P_{t_{r_i}}(p^i_{r_i})\otimes 1 \otimes \cdots \otimes
p_n^i\neq 0$, alors
$\sum_{j=0}^{r_i-2}\mathrm{deg}(p^i_j)+\mathrm{deg}(
p^i_{r_{i}-1}P_{t_{r_i}}(p^i_{r_i}))\geq r_i $, et si
$p_0^i \otimes
\cdots \otimes p^i_{r_{i}-1}I'_{t_{r_i}}(p^i_{r_i})\otimes
x_{t_{r_i}} \otimes \cdots \otimes p_n^i)\neq 0$, alors
$\sum_{j=0}^{r_i-2}\mathrm{deg}(p^i_j)+\mathrm{deg}(
p^i_{r_{i}-1}I'_{t_{r_i}}(p^i_{r_i}))\geq r_i $.

Ceci montre que s'il existe une \'ecriture de $x$ qui satisfait la propri\'et\'e (*M) alors il existe une autre \'ecriture qui satisfait la propri\'et\'e (*M$-1$). Par r\'ecurrence on arrive au fait qu'il existe une \'ecriture de  $x$ qui satisfait la propri\'et\'e (*0), ce qui est \'equivalent \`a dire que $x$ est sup\'erieur. \end{proof}

\section{Une base de morphismes dans le cas g\'en\'eral}
Maintenant qu'on a prouv\'{e} le th\'{e}or\`{e}me \ref{bacan}, on finit en donnant une solution du probl\`eme \ref{pr}, c'est \`a dire, on donne une base explicite de
 $\mathrm{Hom}(\theta_{s_1}\cdots \theta_{s_n},
\theta_{t_1}\cdots \theta_{t_k})$, pour $s_1,\ldots , s_n, t_1, \ldots, t_k \in \mathcal{S}.$.

\begin{defn}
Si $\overline{i}=(i_1, \ldots , i_r) \in I_r$ on d\'{e}finit
$\overline{i}^c=(-i_1+1,\ldots , -i_r+1) \in I_r$. On d\'{e}finit aussi $\overline{i}(\mathrm{op})=(i_r,
\ldots , i_1)$. Finalement on pose $x_{\overline{i}}^g=
x_{t_1}^{i_1} \otimes \cdots \otimes x_{t_r}^{i_r}\otimes 1 \in
\theta_{t_1}\cdots \theta_{t_r}$ et $x_{\overline{i}}^d= 1\otimes
x_{t_1}^{i_1} \otimes \cdots \otimes x_{t_r}^{i_r} \in
\theta_{t_1}\cdots \theta_{t_r}$.
\end{defn}

Soit $\{ f_{\alpha}\}_{\alpha \in A}\subseteq \mathrm{Hom}(\theta_{t_k} \cdots \theta_{t_1}\theta_{s_1}
 \cdots
\theta_{s_n}(-k),R)$ une base de feuilles l\'eg\`eres (noter qu'on dit \og une  \fg \  base parce que pour choisir cette base il faut faire un choix du m\^eme type que dans \ref{para}). Alors avec l'isomorphisme explicite donn\'{e} dans le lemme \ref{bijection}, on trouve
 que
 $$\left\{ m \mapsto \sum_{\overline{i}\in I_k} x_{\overline{i}}^gf_{\alpha}((1\otimes
 x_{\overline{i}^c(\mathrm{op})}^d)\cdot m) \right\}_{\alpha \in A } $$
 est une base de $\mathrm{Hom}(\theta_{s_1} \cdots
\theta_{s_n},\theta_{t_1} \cdots \theta_{t_k})$ comme $R$-module \`a droite.

\section{Application à la catégorie $\mathcal{O}$ de BGG}

 Soit $\mathfrak{g}$ une algèbre de Lie semisimple complexe et $W$ son groupe de Weyl. On rappelle que $\mathcal{O}_0-$proj est la sous-catégorie des objets projectifs dans  le bloc principal $\mathcal{O}_0$ de la catégorie $\mathcal{O}$. Soit $C=R/(R_+^W)$ l'algèbre de coinvariants. On identifie $\mathbb{C}$ à $R/R_+$, et on définit $\mathbf{B}^{\mathbb{C}}$ la  sous-catégorie pleine de $C-$mod, d'objets  les éléments $B\otimes_R \mathbb{C}$, avec $B\in \mathbf{B}$.

Dans l'article \cite{S0}, Soergel construit une  équivalence explicite de catégories :
\begin{equation}
 \mathbb{V} : \mathcal{O}_0-\mathrm{proj}\simeq  \mathbf{B}^{\mathbb{C}}.
\end{equation}

D'autre part, à partir de la proposition 8 de l'article \cite{S1}, on obtient le corollaire suivant :

\begin{cor}
Si $B,B'\in \mathbf{B}$, alors le morphisme canonique
 $$\mathrm{Hom}_{(R,R)}(B,B')\otimes_R \mathbb{C} \simeq \mathrm{Hom}_C(B\otimes_R \mathbb{C}, B'\otimes_R \mathbb{C}) $$
est un isomorphisme.
\end{cor}
Avec ces deux résultats  et la description de $\mathrm{Hom}_{(R,R)}(B,B')$ faite dans la section 6, on obtient explicitement les morphismes dans la catégorie $\mathcal{O}_0$-proj.

On remarque  que $K^b(\mathcal{O}_0-\mathrm{proj})\cong D^b(\mathcal{O}_0-\mathrm{mod}),$ car $\mathrm{gldim}\mathcal{O}_0<\infty$.


\begin{thebibliography}{10}


  \bibitem{BB} \textbf{A.~Beilinson et J.~Bernstein},
 \emph{Localization de $\mathfrak{g}$-modules}, C.R. Acad. Sci. Paris (1) \textbf{292} (1981),
 15-18.

  \bibitem{B}  \textbf{N.~Bourbaki}, \emph{\'El\'ements de
 math\'ematique. Groupes et alg\`ebres de Lie : chapitres 4, 5 et 6}, Hermann, Paris,
 1968.

  \bibitem{BK} \textbf{J.L.~Brylinski et M.~Kashiwara},
 \emph{Kazhdan-Lusztig conjecture and holonomic systems}, Invent. Math. \textbf{64} (1981),
 387-410.

  \bibitem{D} \textbf{M.~Dyer}, \emph{On some generalisations of the
 Kazhdan-Lusztig polynomials for \og universal\fg \  Coxeter systems}, J. Algebra
 \textbf{116 }(1988), no. 2, 353-371.

  \bibitem{F} \textbf{P.~Fiebig}, \emph{The combinatorics of Coxeter
 categories}, preprint {\tt arXiv:math.RT/0512176}, to appear in Trans. Amer. Math. Soc.

\bibitem{GP} \textbf{M.~Geck et G.~Pfeiffer}  \emph{Characters of
 finite Coxeter groups and Iwahori-Hecke algebras}, Oxford University
 Press, Oxford, 2000.
  
\bibitem{H} \textbf{Z.~Haddad}, \emph{A Coxeter group approach to
 Schubert varieties}, Infinite-dimensional groups with applications
 (Berkeley, California
  1984),  Math. Sci. Res. Inst. Publ. \textbf{4}, Springer, New York/Berlin,
 1985, 157-165.

  \bibitem{Hi} \textbf{H.~Hiller}, \emph{Geometry of Coxeter groups},
  Research Notes in Mathematics, \textbf{No. 54}, Pitman, Boston, 1982.

  \bibitem{KL} \textbf{D.~Kazhdan et G.~Lusztig}, \emph{Schubert
 varieties and Poincar\'{e}
duality}, Proc. Symp. Pure Math \textbf{36}, 1980.

  \bibitem{K} \textbf{M.~Khovanov}, \emph{Triply-graded link homology
 and Hochschild homology of Soergel bimodules}, preprint {\tt arXiv:math.GT/0510265}, to appear in  International Journal of Math.

  \bibitem{KR} \textbf{M.~Khovanov et L.~Rozansky}, \emph{Matrix
 factorizations and link homology II}, preprint {\tt arXiv:math.QA/0505056}, to appear in Geometry and Topology.

   
\bibitem{lib} \textbf{N.~Libedinsky} \emph{\'Equivalences entre conjectures de Soergel}, preprint, arXiv:math.RT/0802.3031.

\bibitem{S0} \textbf{W.~Soergel}, \emph{Kategorie  O, perverse Garben, und Moduln über den Koinvarianten zur Weylgruppe}, Journal of the AMS \textbf{3}, (1990), 421-445.

\bibitem{S1} \textbf{W.~Soergel}, \emph{The combinatorics of
 Harish-Chandra bimodules.}, J. Reine
Angew. Math. \textbf{429}, (1992), 49-74.

  \bibitem{S2} \textbf{W.~Soergel}, \emph{On the relation between
 intersection cohomology and representation
theory in positive characteristic}, J. Pure Appl. Algebra \textbf{152}
(2000), no. 1-3, 311-335.

  \bibitem{S3} \textbf{W.~Soergel}, \emph{Kazhdan-Lusztig polynomials
 and indecomposable bimodules over
polynomial rings}, Journal of the Inst. of Math. Jussieu  \textbf{6(3)}, (2007), 501-525.


\end{thebibliography}
\end{document}